\documentclass[onecolumn]{autart}

\usepackage{graphicx}

\begin{document}
\input{amssym}
\begin{frontmatter}

\title{Group Analysis of Born-Infeld Equation}
%\thanksref{footnoteinfo}} % Title, preferably not more
                                                % than 10 words.

%\thanks[footnoteinfo]{This paper was not presented at any IFAC
%meeting. Corresponding author M.~T.~Cicero. Tel. +XXXIX-VI-mmmxxi.

%Fax +XXXIX-VI-mmmxxv.}

\author[]{Mehdi Nadjafikhah}\ead{m\_nadjafikhah@iust.ac.ir},    % Add the
\author[MN]{Seyed Reza Hejazi}\ead{reza\_hejazi@iust.ac.ir}       % e-mail address
%\author[Baiae]{Publius Maro Vergilius}\ead{vergilius@culture.ir}  % (ead) as shown

\address[MN]{Iran University of Science and Technology, Narmak, Tehran}  % Please supply
%\address[RH]{Iran University of Science and Technology, Narmak, Tehran}             % full addresses
%\address[Baiae]{The White House, Baiae}        % here.

                             % keyword list or with the
                                          % help of the Automatica
                                          % keyword wizard

\begin{abstract}                          % Abstract of not more than 200 words.
Lie symmetry group method is applied to study the Born-Infeld
equation. The symmetry group and its optimal system are given, and
group invariant solutions associated to the symmetries are
obtained. Finally the structure of the Lie algebra symmetries is
determined.
\end{abstract}

\end{frontmatter}

\section{Introduction}
$~~~~~$The method of point transformations are a powerful tool in
order to find exact solutions for nonlinear partial differential
equations. It happens that many PDE's of physical importance are
nonlinear and Lie classical symmetries admitted by nonlinear PDE's
are useful for finding invariant solutions.\\
$~~~~~$In physics, the Born-Infeld theory is a nonlinear
generalization of electromagnetism \cite{[6],[7]}. The model is
named after physicists Max Born (1882-1970) and Leopold Infeld
(1898-1968) who first proposed it. In physics, it is a particular
example of what is usually known as a nonlinear electrodynamics.
It was historically introduced in the 30's to remove the
divergence of the electron's self-energy in classical
electrodynamics by introducing an upper bound of the electric
field at the origin. The Born-Infeld electrodynamics possesses a
whole series of physically interesting properties: First of all
the total energy of the electromagnetic field is finite and the
electric field is regular everywhere. Second it displays good
physical properties concerning wave propagation, such as the
absence of shock waves and birefringence. A field theory showing
this property is usually called completely exceptional and
Born-Infeld theory is the only completely exceptional regular
nonlinear electrodynamics. Finally (and more technically)
Born-Infeld theory can be seen as a covariant generalization of
Mie's theory, and very close to Einstein's idea of introducing a
nonsymmetric metric tensor with the symmetric part corresponding
to the usual metric tensor and the antisymmetric to the
electromagnetic field tensor. During the 1990s there was a revival
of interest on Born-Infeld theory and its nonabelian extensions as
they were found in some limits of string theory.
\section{Lie Symmetries of the Equation}
A PDE with $p-$independent and $q-$dependent variables has a Lie
point transformations
\begin{eqnarray*}
\widetilde{x}_i=x_i+\varepsilon\xi_i(x,u)+{\mathcal
O}(\varepsilon^2),\qquad
\widetilde{u}_{\alpha}=u_\alpha+\varepsilon\varphi_\alpha(x,u)+{\mathcal
O}(\varepsilon^2)
\end{eqnarray*}
where
$\displaystyle{\xi_i=\frac{\partial\widetilde{x}_i}{\partial\varepsilon}\Big|_{\varepsilon=0}}$
for $i=1,...,p$ and
$\displaystyle{\varphi_\alpha=\frac{\partial\widetilde{u}_\alpha}{\partial\varepsilon}\Big|_{\varepsilon=0}}$
for $\alpha=1,...,q$. The action of the Lie group can be
considered by its associated infinitesimal generator
\begin{eqnarray}\label{eq:18}
\textbf{v}=\sum_{i=1}^p\xi_i(x,u)\frac{\partial}{\partial{x_i}}+\sum_{\alpha=1}^q\varphi_\alpha(x,u)\frac{\partial}{\partial{u_\alpha}}
\end{eqnarray}
on the total space of PDE (the space containing independent and
dependent variables). Furthermore, the characteristic of the
vector field (\ref{eq:18}) is given by
\begin{eqnarray*}
Q^\alpha(x,u^{(1)})=\varphi_\alpha(x,u)-\sum_{i=1}^p\xi_i(x,u)\frac{\partial
u^\alpha}{\partial x_i},
\end{eqnarray*}
and its $n-$th prolongation is determined by
\begin{eqnarray*}
\textbf{v}^{(n)}=\sum_{i=1}^p\xi_i(x,u)\frac{\partial}{\partial
x_i}+\sum_{\alpha=1}^q\sum_{\sharp
J=j=0}^n\varphi^J_\alpha(x,u^{(j)})\frac{\partial}{\partial
u^\alpha_J},
\end{eqnarray*}
where
$\varphi^J_\alpha=D_JQ^\alpha+\sum_{i=1}^p\xi_iu^\alpha_{J,i}$.
($D_J$ is the total derivative operator describes in
(\ref{eq:19})).

The aim is to analysis the Lie point symmetry structure of the
Born-Infeld equation, which is
\begin{equation}\label{eq:1}
(1-u_t^2)u_{xx}+2u_xu_tu_{xt}-(1+u_x^2)u_{tt}=0,
\end{equation}
where $u$ is a smooth function of $\displaystyle{(x,t)}$. Let us
consider a one-parameter Lie group of infinitesimal
transformations $(x,t,u)$ given by
\begin{eqnarray*}
\widetilde{x}=x+\varepsilon\xi_1(x,t,u)+{\mathcal
O}(\varepsilon^2),\qquad
\widetilde{t}=t+\varepsilon\xi_2(x,t,u)+{\mathcal
O}(\varepsilon^2),\qquad
\widetilde{u}=u+\varepsilon\eta(x,t,u)+{\mathcal
O}(\varepsilon^2),
\end{eqnarray*}
where $\varepsilon$ is the group parameter. Then one requires that
this transformations leaves invariant the set of solutions of the
Eq. (\ref{eq:1}). This yields to the linear system of equations
for the infinitesimals $\xi_1(x,t,u)$, $\xi_2(x,t,u)$ and
$\eta(x,t,u)$. The Lie algebra of infinitesimal symmetries is the
set of vector fields in the form of
\begin{eqnarray*}
\textbf{v}=\xi_1(x,t,u){\partial_x}+\xi_2(x,t,u){\partial_t}+\eta(x,t,u){\partial_u}.
\end{eqnarray*}
This vector field has the second prolongation
\begin{eqnarray*}
\textbf{v}^{(2)}=\textbf{v}+\varphi^x\partial_{u_x}+\varphi^t\partial{u_t}+\varphi^{xx}\partial_{u_{xx}}+\varphi^{xt}\partial_{u_{xx}}+\varphi^{tt}\partial_{u_{tt}}
\end{eqnarray*}
with the coefficients
\begin{eqnarray*}
\varphi^x&=&D_x(\varphi-\xi_1u_x-\xi_2u_t)+\xi_1u_{xx}+\xi_2u_{xt},\\
\varphi^t&=&D_t(\varphi-\xi_1u_x-\xi_2u_t)+\xi_1u_{xt}+\xi_2u_{tt},\\
\varphi^{xx}&=&D^2_x(\varphi-\xi_1u_x-\xi_2u_t)+\xi_1u_{xxx}+\xi_2u_{xxt},\\
\varphi^{xt}&=&D_xD_t(\varphi-\xi_1u_x-\xi_2u_t)+\xi_1u_{xxt}+\xi_2u_{xtt},\\
\varphi^{tt}&=&D^2_t(\varphi-\xi_1u_x-\xi_2u_t)+\xi_1u_{xtt}+\xi_2u_{ttt},\\
\end{eqnarray*}
where the operators $D_x$ and $D_t$ denote the total derivative
with respect to $x$ and $t$:
\begin{eqnarray}\label{eq:19}
D_x&=&\partial_x+u_x\partial_u+u_{xx}\partial_{u_x}+u_{xt}\partial_{u_t}+\cdots,\\
D_t&=&\partial_t+u_t\partial_u+u_{tt}\partial_{u_t}+u_{xt}\partial_{u_x}+\cdots.\nonumber
\end{eqnarray}
Using the invariance condition, i.e., applying the second
prolongation $\textbf{v}^{(2)}$ to Eq. (\ref{eq:1}), the following
system of 10 determining equations yields:
\begin{eqnarray*}
\begin{array}{lclclclclc}
{\xi_2}_{xx}=0,&&{\xi_2}_{xu}=0,&&{\xi_2}_{tt}=0,&&{\xi_2}_{uu}=0,&&{\xi_1}_x={\xi_2}_t,\\
{\xi_1}_t={\xi_2}_x,&&{\xi_1}_u=-\eta_x,&&{\xi_2}_t=\eta_u,&&{\xi_2}_u=\eta_t,&&{\xi_2}_{tu}=-\eta_{xx}.
\end{array}
\end{eqnarray*}
The solution of the above system gives the following coefficients
of the vector field $\textbf{v}$:
\begin{eqnarray*}
\xi_1=c_1+c_4t-c_5u+c_7x,\qquad \xi_2=c_2+c_4x+c_6u+c_7t,\qquad
\eta=c_3+c_5x+c_6t+c_7t,
\end{eqnarray*}
where $c_1,...,c_7$ are arbitrary constants, thus the Lie algebra
${\goth g}$ of the Born-Infeld equation is spanned by the seven
vector fields
\begin{eqnarray*}
\begin{array}{lclclclc}
\textbf{v}_1=\partial_x,&& \textbf{v}_2=\partial_t,&&
\textbf{v}_3=\partial_u,&&\textbf{v}_4=t\partial_x+x\partial_t,\\
\textbf{v}_5=-u\partial_x+x\partial_u,&&\textbf{v}_6=u\partial_t+t\partial_u,&&
\textbf{v}_7=x\partial_x+t\partial_t+u\partial_u,
\end{array}
\end{eqnarray*}
which $\textbf{v}_1,\textbf{v}_2$ and $\textbf{v}_3$ are
translation on $x,t$ and $u$, $\textbf{v}_5$ is rotation on $u$
and $x$ and $\textbf{v}_7$ is scaling on $x,t$ and $u$. The
commutation relations between these vector fields is given by the
table 1, where entry in row $i$ and column $j$ representing
$[\textbf{v}_i,\textbf{v}_j]$.
\begin{table}
\caption{Commutation relations of $\goth g$ }\label{table:1}
\vspace{-0.3cm}\begin{eqnarray*}\hspace{-0.75cm}\begin{array}{llllllll}
\hline
  [\,,\,]       &\hspace{1.1cm}\textbf{v}_1  &\hspace{0.5cm}\textbf{v}_2  &\hspace{0.5cm}\textbf{v}_3 &\hspace{0.5cm}\textbf{v}_4 &\hspace{0.5cm}\textbf{v}_5 &\hspace{0.5cm}\textbf{v}_6&\hspace{0.5cm}\textbf{v}_7  \\ \hline
  \textbf{v}_1  &\hspace{1.1cm} 0            &\hspace{0.5cm} 0            &\hspace{0.5cm}0            &\hspace{0.5cm}\textbf{v}_3 &\hspace{0.5cm}\textbf{v}_3 &\hspace{0.5cm}0           &\hspace{0.5cm}\textbf{v}_1     \\
  \textbf{v}_2  &\hspace{1.1cm} 0            &\hspace{0.5cm} 0            &\hspace{0.5cm}0            &\hspace{0.5cm}\textbf{v}_1 &\hspace{0.5cm}0            &\hspace{0.5cm}\textbf{v}_3&\hspace{0.5cm}\textbf{v}_2     \\
  \textbf{v}_3  &\hspace{1.1cm} 0            &\hspace{0.5cm} 0            &\hspace{0.5cm}0            &\hspace{0.5cm}0            &\hspace{0.3cm}-\textbf{v}_1&\hspace{0.5cm}\textbf{v}_2&\hspace{0.5cm}\textbf{v}_3     \\
  \textbf{v}_4  &\hspace{0.8cm}-\textbf{v}_2 &\hspace{0.3cm} -\textbf{v}_1&\hspace{0.5cm}0            &\hspace{0.5cm}0            &\hspace{0.5cm}\textbf{v}_6 &\hspace{0.5cm}\textbf{v}_5&\hspace{0.5cm}0     \\
  \textbf{v}_5  &\hspace{0.8cm}-\textbf{v}_3 &\hspace{0.5cm} 0            &\hspace{0.5cm}\textbf{v}_1 &\hspace{0.3cm}-\textbf{v}_6&\hspace{0.5cm}0            &\hspace{0.5cm}\textbf{v}_4&\hspace{0.5cm}0     \\
  \textbf{v}_6  &\hspace{1.1cm} 0            &\hspace{0.3cm} -\textbf{v}_3&\hspace{0.3cm}-\textbf{v}_2&\hspace{0.3cm}-\textbf{v}_5&\hspace{0.3cm}-\textbf{v}_4&\hspace{0.5cm}0           &\hspace{0.5cm}0     \\
  \textbf{v}_7  &\hspace{0.8cm}-\textbf{v}_1&\hspace{0.3cm} -\textbf{v}_2&\hspace{0.3cm}-\textbf{v}_3&\hspace{0.5cm}0            &\hspace{0.5cm}0            &\hspace{0.5cm}0           &\hspace{0.5cm}0     \\
  \hline\end{array}\end{eqnarray*}\end{table}

The one-parameter groups $G_i$ generated by the base of $\goth g$
are given in the following table.
\begin{eqnarray*}
\begin{array}{lclc}
G_1:(x+\varepsilon,t,u),&& G_2:(x,t+\varepsilon,u),\\
G_3:(x,t,u+\varepsilon),&&
G_4:\Big(x\cosh\varepsilon+t\sinh\varepsilon,x\sinh\varepsilon+t\cosh\varepsilon,u\Big)\\[2mm]
G_5:(-u\sin\varepsilon+x\cos\varepsilon,t,x\sin\varepsilon+u\cos\varepsilon),&&
G_6:\Big(x,t\cosh\varepsilon+u\sinh\varepsilon,t\sinh\varepsilon+u\cosh\varepsilon\Big),\\[2mm]
G_7:(xe^{\varepsilon},te^{\varepsilon},ue^{\varepsilon}).
\end{array}
\end{eqnarray*}
Since each group $G_i$ is a symmetry group and if $u=f(x,t)$ is a
solution of the Born-Infeld equation, so are the functions
\begin{eqnarray*}
\begin{array}{lll}
u_1=f(x+\varepsilon,t),\qquad u_2=f(x,t+\varepsilon),&&
u_3=f(x,t)-\varepsilon,\\
u_4=f\Big(x\cosh\varepsilon-t\sinh\varepsilon,-x\sinh\varepsilon+t\cosh\varepsilon\Big),&&
u_5=\sec\varepsilon f(u\sin\varepsilon+x\cos\varepsilon,t)+x\sin\varepsilon,\\
u_6=\sec\mbox{h}\varepsilon
f(x,t\cosh\varepsilon-u\sinh\varepsilon)+t\sinh\varepsilon ,&&
u_7=e^{-\varepsilon}f(e^{-\varepsilon}x,e^{-\varepsilon}t).
\end{array}
\end{eqnarray*}
where $\varepsilon$ is a real number. Here we can find the general
group of the symmetries by considering a general linear
combination $c_1\textbf{v}_1+\cdots+c_1\textbf{v}_6$ of the given
vector fields. In particular if $g$ is the action of the symmetry
group near the identity, it can be represented in the form
$g=\exp(\varepsilon_7\textbf{v}_7)\cdots\exp(\varepsilon_1\textbf{v}_1)$.
\section{Symmetry reduction for Born-Infeld equation}
$~~~~~$The first advantage of symmetry group method is to
construct new solutions from known solutions. Neither the first
advantage nor the second will be investigated here, but symmetry
group method will be applied to the Eq. (\ref{eq:1}) to be
connected directly to some order differential equations. To do
this, a particular linear combinations of infinitesimals are
considered and their corresponding invariants are determined. The
Born-Infeld equation expressed in the coordinates $(x,t)$, so to
reduce this equation is to search for its form in specific
coordinates. Those coordinates will be constructed by searching
for independent invariants $(y,v)$ corresponding to an
infinitesimal generator. so using the chain rule, the expression
of the equation in the new coordinate allows us to the reduced
equation. Here we will obtain some invariant solutions with
respect to symmetries. First we obtain the similarity variables
for each term of the Lie algebra $\goth g$, then we use this
method to reduced the PDE and find the invariant solutions. All
results are coming in the following table.
\begin{eqnarray*}
\begin{array}{lclclcl}
\mbox{\textbf{vector field}}&&\mbox{\textbf{invariant
function}}&&\mbox{\textbf{invariant
transformations}}&&\mbox{\textbf{invariant solution}}\\
\textbf{v}_1&&\Phi(t,u)&&v(y)=u(x,t),y=t&&u=c_1t+c_2\\
\textbf{v}_2&&\Phi(x,u)&&v(y)=u(x,t),y=x&&u=c_1x+c_2\\
\textbf{v}_3&&\Phi(x,t)&&\cdots&&\mbox{translation of all solution
is the invariant}\\
&&&&&&\mbox{solution}\\
\textbf{v}_4&&\Phi(-x^2+t^2,u)&&v(y)=u(x,t),y=-x^2+t^2&&u=\pm
c_1\arctan\Big(\frac{x^2-t^2+2c_1^2}{\sqrt{(x^2-t^2)(-x^2+t^2-4c_1^2)}}\Big)+c_2,\\
\textbf{v}_5&&\Phi(t,x^2+u^2)&&v(y)=x^2+u(x,t)^2,y=t&&u=\pm\frac{1}{2}\sqrt{2c_1e^{\frac{t+c_2}{c_1}}+8c_1^3e^{-\frac{t+c_2}{c_1}}-4x^2-8c_1^2}\\
&&&&&&u=\pm\frac{1}{2}\sqrt{2c_1e^{-\frac{t+c_2}{c_1}}+8c_1^3e^{\frac{t+c_2}{c_1}}-4x^2-8c_1^2}\\
\textbf{v}_6&&\Phi(x,-t^2+u^2)&&v(y)=-t^2+u(x,t)^2,y=x&&u=\pm\frac{1}{2}\sqrt{2c_1e^{\frac{x+c_2}{c_1}}+8c_1^3e^{-\frac{x+c_2}{c_1}}+4x^2+8c_1^2}\\
&&&&&&u=\pm\frac{1}{2}\sqrt{2c_1e^{-\frac{x+c_2}{c_1}}+8c_1^3e^{-\frac{x+c_2}{c_1}}+4x^2+8c_1^2}\\
\textbf{v}_7&&\displaystyle{\Phi\Big(\frac{t}{x},\frac{u(x,t)}{x}\Big)}&&\displaystyle{v(y)=\frac{u(x,t)}{x}},\displaystyle{y=\frac{t}{x}}&&u=c_1x+c_2t,u=\pm\sqrt{-x^2+t^2}
\end{array}
\end{eqnarray*}
\section{Optimal system of Born-Infeld equation}
$~~~~~$As is well known, the theoretical Lie group method plays an
important role in finding exact solutions and performing symmetry
reductions of differential equations. Since any linear combination
of infinitesimal generators is also an infinitesimal generator,
there are always infinitely many different symmetry subgroups for
the differential equation. So, a mean of determining which
subgroups would give essentially different types of solutions is
necessary and significant for a complete understanding of the
invariant solutions. As any transformation in the full symmetry
group maps a solution to another solution, it is sufficient to
find invariant solutions which are not related by transformations
in the full symmetry group, this has led to the concept of an
optimal system \cite{[3]}. The problem of finding an optimal
system of subgroups is equivalent to that of finding an optimal
system of subalgebras. For one-dimensional subalgebras, this
classification problem is essentially the same as the problem of
classifying the orbits of the adjoint representation. This problem
is attacked by the naive approach of taking a general element in
the Lie algebra and subjecting it to various adjoint
transformations so as to simplify it as much as possible. The idea
of using the adjoint representation for classifying
group-invariant solutions is due to \cite{[4]} and \cite{[5]}.

The adjoint action is given by the Lie series
\begin{eqnarray}\label{eq:9}
\mbox{Ad}(\exp(\varepsilon\textbf{v}_i)\textbf{v}_j)=\textbf{v}_j-\varepsilon[\textbf{v}_i,\textbf{v}_j]+\frac{\varepsilon^2}{2}[\textbf{v}_i,[\textbf{v}_i,\textbf{v}_j]]-\cdots,
\end{eqnarray}
where $[\textbf{v}_i,\textbf{v}_j]$ is the commutator for the Lie
algebra, $t$ is a parameter, and $i,j=1,\cdots,10$. Let
$F^{\varepsilon}_i:{\goth g}\rightarrow{\goth g}$ defined by
$\textbf{v}\mapsto\mbox{Ad}(\exp(\varepsilon\textbf{v}_i)\textbf{v})$
is a linear map, for $i=1,\cdots,7$. The matrices
$M^\varepsilon_i$ of $F^\varepsilon_i$, $i=1,\cdots,7$, with
respect to basis $\{\textbf{v}_1,\cdots,\textbf{v}_7\}$ are
\begin{eqnarray*}
\begin{array}{lcl}
M^\varepsilon_1=\small\left(\begin{array}{ccccccc}
1&0&0&0&0&0&0\\0&1&0&0&0&0&0\\0&0&1&0&0&0&0\\0
&-\varepsilon&0&1&0&0&0\\0&0&-\varepsilon&0&1&0&0\\0&0&0&0&0&1&0\\-\varepsilon&0&0&0&0&0&1\end{array}
\right),&& M^\varepsilon_2=\small\left(\begin{array}{ccccccc}
1&0&0&0&0&0&0\\0&1&0&0&0&0&0\\0&0&1&0&0&0&0\\-\varepsilon
&0&0&1&0&0&0\\0&0&0&0&1&0&0\\0&0&-\varepsilon&0&0&1&0\\0&-\varepsilon&0&0&0&0&1\end{array}
\right),\\ M^\varepsilon_3=\small\left(\begin{array}{ccccccc}
1&0&0&0&0&0&0\\0&1&0&0&0&0&0\\0&0&1&0&0&0&0\\0 &0&0&1&0&0&0\\
\varepsilon&0&0&0&1&0&0\\0&-\varepsilon&0&0&0&1&0\\0&0&-\varepsilon&0&0&0&1\end{array}
\right),&& M^\varepsilon_4=\small\left(\begin{array}{ccccccc}
\cosh\varepsilon&\sinh\varepsilon&0&0&0&0&0\\\sinh\varepsilon&\cosh\varepsilon&0&0&0&0&0\\0&0&1&0&0&0&0\\0&0&0&1&0&0&0\\0&0&0&0&\cosh\varepsilon&-\sinh\varepsilon&0
\\0&0&0&0&-\sinh\varepsilon&\cosh\varepsilon&0\\0&0&0&0&0&0&1\end{array}
\right),\\ M^\varepsilon_5=\small\left(\begin{array}{ccccccc}
\cos\varepsilon&0&\sin\varepsilon&0&0&0&0\\0&1&0&0&0&0&0\\-\sin\varepsilon&0&\cos\varepsilon&0&0&0&0\\0&0&0&\cos\varepsilon&0&\sin\varepsilon&0\\0&0&0&0&1&0&0\\0&0&0&-\sin\varepsilon&0&\cos\varepsilon&0\\0&0&0&0&0&0&1\end{array}
\right),&& M^\varepsilon_6=\small\left(\begin{array}{ccccccc}
1&0&0&0&0&0&0\\0&\cosh\varepsilon&\sin\varepsilon&0&0&0&0\\0&\sin\varepsilon&\cos\varepsilon&0&0&0&0\\0
&0&0&\cosh\varepsilon&\sin\varepsilon&0&0\\0&0&0&\sinh\varepsilon&\cosh\varepsilon&0&0\\0&0&0&0&0&1&0\\0&0&0&0&0&0&1\end{array}
\right),\\ M^\varepsilon_7=\small\left(\begin{array}{ccccccc}
e^\varepsilon&0&0&0&0&0&0\\0&e^\varepsilon&0&0&0&0&0\\0&0&e^\varepsilon&0&0&0&0\\0
&0&0&1&0&0&0\\0&0&0&0&1&0&0\\0&0&0&0&0&1&0\\0&0&0&0&0&0&1\end{array}
\right),
\end{array}
\end{eqnarray*}
by acting these matrices on a vector field $\textbf{v}$
alternatively we can  show that a one-dimensional optimal system
of ${\goth g}$ is given by
\begin{eqnarray*}
\begin{array}{lll}
X_1=\textbf{v}_1,&& X_2=\textbf{v}_3,\\
X_3=a(\textbf{v}_1-\textbf{v}_7)+\textbf{v}_2,&& X_4=a\textbf{v}_2-b\textbf{v}_3+c\textbf{v}_7,\\
X_5=\textbf{v}_1+a\textbf{v}_2+\textbf{v}_3-\textbf{v}_6,&& X_6=\textbf{v}_1+a\textbf{v}_2+b(\textbf{v}_3-\textbf{v}_6),\\
X_7=a\textbf{v}_1+b\textbf{v}_2-c\textbf{v}_3-d\textbf{v}_6,&&
X_8=a(\textbf{v}_1-\textbf{v}_4)+b(\textbf{v}_2-\textbf{v}_4)+\textbf{v}_3+c\textbf{v}_4.
\end{array}
\end{eqnarray*}
In the next section we will find the invariant solutions with
respect to the symmetries and optimal system.
\section{Lie Algebra Structure}
$~~~~~$In this part, we determine the structure of symmetry Lie
algebra of the Born-Infeld equation.\\
The Lie algebra $\goth g$ is not solvable and semisimple, because
if $\displaystyle{{\goth g}^{(1)}=\mbox{Span}_{\Bbb
R}\{\textbf{v}_i,[\textbf{v}_i,\textbf{v}_j]\}_{i,j}}$ be the
derived of ${\goth g}$ we have
\begin{eqnarray*}
{\goth g}^{(1)}=\mbox{Span}_{\Bbb
R}\{\textbf{v}_1,...,\textbf{v}_7\} ={\goth g},
\end{eqnarray*}
but it has a \textit{Levi decomposition} in the form of
\begin{eqnarray}\label{eq:16}
{\goth g}={\goth r}\ltimes{\goth g}_1,
\end{eqnarray}
where $\displaystyle{{\goth
r}=\mbox{Span}_{\Bbb{R}}\{\textbf{v}_1,\textbf{v}_2,\textbf{v}_3,\textbf{v}_7\}}$
is the radical (the largest solvable ideal) of ${\goth g}$, and
$\displaystyle{{\goth g}_1=\mbox{Span}_{\Bbb
R}\{\textbf{v}_4,\textbf{v}_5,\textbf{v}_6\}}$ is the semisimple
and nonsolvable subalgebra of ${\goth g}$. So the Levi
decomposition of symmetry Lie algebra for Born-Infeld equation
gives the quotient structure
\begin{eqnarray}\label{eq:17}
{\overline{\goth g}}={\goth g}/{\goth r}.
\end{eqnarray}
If $\textbf{w}_i=\textbf{v}_i+{\goth r}$ are the members of
quotient algebra, the commutators table for ${\overline{\goth g}}$
are given in table 2.
\begin{table}
\caption{Commutation relations of $\goth g$ }\label{table:2}
\vspace{-0.3cm}\begin{eqnarray*}\hspace{-0.75cm}\begin{array}{llll}
\hline
  [\,,\,]       &\hspace{1.1cm}\textbf{w}_1  &\hspace{0.5cm}\textbf{w}_2  &\hspace{0.5cm}\textbf{w}_3  \\ \hline
  \textbf{w}_1  &\hspace{1.4cm} 0            &\hspace{0.7cm} \textbf{w}_3            &\hspace{0.5cm}\textbf{w}_2    \\
  \textbf{w}_2  &\hspace{1.1cm} -\textbf{w}_3            &\hspace{0.8cm} 0            &\hspace{0.5cm}\textbf{w}_1  \\
  \textbf{w}_3  &\hspace{1.1cm} -\textbf{w}_2            &\hspace{0.5cm} -\textbf{w}_1            &\hspace{0.6cm}0    \\
  \hline\end{array}\end{eqnarray*}\end{table}

$~~~~~$Finally, we have some analysis on the structure of
(\ref{eq:16}) and (\ref{eq:17}) with some important objects in
algebra. We know that he centralizer of a set of vectors $\goth g$
in a subalgebra ${\goth h}$ is the subalgebra of vectors in
${\goth h}$ which commute with all the vectors in $\goth g$. With
attentive to (\ref{eq:16}), $\goth r$ has no any nontrivial
centralizer and it is the only minimal ideal containing itself.
But this is not true for ${\goth g}_1$, because its centralizer
has a member which is $\textbf{v}_7$, and the minimal ideal
containing ${\goth g}_1$ is spanned by
$\{\textbf{v}_1,...,\textbf{v}_6\}$.
\section{Conclusion}
In this article group classification of Born-Infeld equation and
the algebraic structure of the symmetry group is considered.
Classification of one-dimensional subalgebra is determined by
constructing one-dimensional optimal system. Some invariant
solutions are fined in the sequel and the Lie algebra structure of
symmetries is found.
%

%\bibliographystyle{plain}        % Include this if you use bibtex
%\bibliography{autosam}           % and a bib file to produce the
                                 % bibliography (preferred). The
                                 % correct style is generated by
                                 % Elsevier at the time of printing.

%\begin{thebibliography}{99}     % Otherwise use the
                                 % thebibliography environment.
                                 % Insert the full references here.
                                 % See a recent issue of Automatica
                                 % for the style.
%  \bibitem[Heritage, 1992]{Heritage:92}
%     (1992) {\it The American Heritage.
%     Dictionary of the American Language.}
%     Houghton Mifflin Company.
%  \bibitem[Able, 1956]{Abl:56}
%     B.~C.~Able (1956). Nucleic acid content of macroscope.
%     {\it Nature 2}, 7--9.
%  \bibitem[Able {\em et al.}, 1954]{AbTaRu:54}
%     B.~C. Able, R.~A. Tagg, and M.~Rush (1954).
%     Enzyme-catalyzed cellular transanimations.
%     In A.~F.~Round, editor,
%     {\it Advances in Enzymology Vol. 2} (125--247).
%     New York, Academic Press.
%  \bibitem[R.~Keohane, 1958]{Keo:58}
%     R.~Keohane (1958).
%     {\it Power and Interdependence:
%     World Politics in Transition.}
%     Boston, Little, Brown \& Co.
%  \bibitem[Powers, 1985]{Pow:85}
%     T.~Powers (1985).
%     Is there a way out?
%     {\it Harpers, June 1985}, 35--47.

%\end{thebibliography}

%\appendix
%\section{A summary of Latin grammar}    % Each appendix must have a short title.
%\section{Some Latin vocabulary}         % Sections and subsections are supported
                                        % in the appendices.
\end{document}